\documentclass[reqno,12pt]{amsart}

\usepackage{amsmath}
\usepackage{amsfonts}
\usepackage{amssymb}
\usepackage{eufrak}
\usepackage{amscd}
\usepackage{amsthm}
\usepackage{amstext}
\usepackage[all]{xy}

 \newcommand{\K}{{\mathbb K}}
 
   \newcommand{\Sp}{\operatorname{Sp}}

 \newcommand{\Ext}{\operatorname{Ext}}

\newcommand{\id}{\operatorname{id}}

   \theoremstyle{plain}
   \newtheorem{thm}{Theorem}
   \newtheorem{prop}[thm]{Proposition}
   
   \newtheorem{cor}[thm]{Corollary}
   \theoremstyle{definition}
   
   \newtheorem{dfn}[thm]{Definition}
   
   \newtheorem{conj}[thm]{Conjecture}
   \theoremstyle{remark}

 \title[]%
{Asymptotic representations of the reduced $C^*$-algebra of a free
group: an example}
\author{V. Manuilov}
\address{Department of Mechanics and Mathematics, Moscow State
University, Leninskie Gory, Moscow, 119992, Russia {\sl and\ }
Harbin Institute of Technology, Harbin, P. R. China}
\email{manuilov@mech.math.msu.su}
\thanks{The author acknowledges partial support from
RFFI, grant No. 05-01-00923}

\begin{document}

\begin{abstract}
We give an example of a non-trivial asymptotic representation of
the reduced $C^*$-algebra of a free group. This example allows to
evaluate the asymptotic tensor $C ^*$-norm of some elements in
tensor product $C^*$-algebras and to show semi-invertibility of
the non-invertible extension of $C^*_r(\mathbb F_2)$ considered by
Haagerup and Thorbj\o rnsen.

\end{abstract}

\maketitle

\section{Introduction}

Asymptotic homomorphisms of $C^*$-algebras were first defined and
studied in \cite{C-H} in relation to topological properties of
$C^*$-algebras. The most important and the best known case is the
case of asymptotic homomorphisms from a suspended $C^*$-algebra
$SA$ to the $C^*$-algebra $\K$ of compact operators, since the
homotopy classes of those are the $K$-homology of $A$, the
$E$-theory. Asymptotic homomorphisms to other $C^*$-algebras are
less known. For example, it is known that any asymptotic
homomorphism to the Calkin algebra is homotopic to a genuine
homomorphism \cite{M-Crelle}. Even less is known about asymptotic
homomorphisms to $\mathbb B(H)$, where there is no topological
obstruction (recall that the $K$-groups of $\mathbb B(H)$ are
trivial). Such asymptotic homomorphisms are called {\it asymptotic
representations} and were first studied in relation to the
asymptotic tensor product $C^*$-algebras \cite{MT6} and to
semi-invertibility of $C^*$-algebra extensions \cite{MT7}.

The aim of this paper is to give an example of an asymptotic
representation of the reduced group $C^*$-algebra of the free
group $\mathbb F_2$ on two generators. This was made possible due
to two interesting families of representations. The first family
was constructed by M.~Pimsner and D.~Voiculescu in \cite{PV}, cf.
next section, and the second family (sequence) was discovered by
U.~Haagerup and S.~Thorbj\o rnsen in \cite{Haagerup-T} by free
probability methods. We combine these two families into an
asymptotic representation of $C^*_r(\mathbb F_2)$.

U.~Haagerup and S.~Thorbj\o rnsen used their sequence of
representations \cite{Haagerup-T} to construct an extension of
$C^*_r(\mathbb F_2)$ which is not invertible in the BDF functor
$\Ext(C^*_r(\mathbb F_2))$ of L. Brown, R. Douglas and P. Fillmore
\cite{BDF}. We are able to show that their extension is
semi-invertible \cite{MT3}, and represents the trivial element in
the modified $\Ext$ functor. Another application of our asymptotic
representation is evaluation of asymptotic tensor $C^*$-norm on
tensor products of $C^*_r(\mathbb F_2)$ with some other
$C^*$-algebras related to the free group. We also show that such
asymptotic tensor products behave differently from the case of
groups with the property T of Kazhdan.

\section{A family of representations of $\mathbb F_2$}

For a free group $\mathbb F_2$, M.~Pimsner and D.~Voiculescu
\cite{PV} have constructed a family $(\pi_s)_{s\in[0,1]}$ of
representations on $l^2(\mathbb F_2)$ with the following
properties:
 \begin{enumerate}
 \item
$\pi_0=\lambda$ is the regular representation;
 \item
$\pi_1=\theta\oplus\lambda^{(2)}$, where $\theta$ is the
one-dimensional trivial representation and
$\lambda^{(2)}=\lambda\oplus\lambda$;
 \item
the map $s\mapsto\pi_s(g)$ is norm-continuous for any $g\in\mathbb
F_2$;
 \item
$\pi_s(g)-\pi_0(g)$ is compact for any $g\in\mathbb F_2$ and any
$s\in[0,1]$.
 \end{enumerate}
Importance of this homotopy, which connects the regular
representation with a representation containing the trivial
representation as a direct summand, was emphasized by E.~C.~Lance
\cite{Lance}. A natural specific example of such a family $\pi_s$
of representations of $\mathbb F_2$ with additional properties is
given in \cite{Pytlik-Szwarc}.


For a representation $\pi$ of $\mathbb F_2$ let us denote by
$C^*_\pi(\mathbb F_2)$ the $C^*$-algebra generated by all
$\pi(g)$, $g\in\mathbb F_2$. We write $\pi(a)\in C^*_\pi(\mathbb
F_2)$ for $a\in C^*(\mathbb F_2)$ and make no distinction between
$a$ and its image in $C^*_\pi(\mathbb F_2)$ when there is no
confusion.



 \begin{prop}\label{exact}
$C^*_{\pi_s}(\mathbb F_2)$ is an exact $C^*$-algebra for any
$s\in[0,1]$.
 \end{prop}
 \begin{proof}
For $s=0$ it is known that $C^*_r(\mathbb F_2)$ is exact, so let
us consider the case $s\in(0,1]$. It follows from the property (4)
of the family $(\pi_s)$ that
 $$
C^*_{\pi_s}(\mathbb F_2)\subset C^*_r(\mathbb F_2)+\K
 $$
for any $s\in(0,1]$. The extension
 $$
0\to\K\to C^*_r(\mathbb F_2)+\K\to C^*_r(\mathbb F_2)\to 0
 $$
is obviously split. So, by results of E. Kirchberg (cf.
\cite{Kirchberg}), $C^*_r(\mathbb F_2)+\K$ is exact, hence its
$C^*$-subalgebras are exact too.

 \end{proof}

\section{Example of an asymptotic representation of $C^*_r(\mathbb F_2)$}

Recall that, by the Fell's absorbtion principle \cite{Fell}, for
any representation $\sigma$ of $\mathbb F_2$, one has
$\lambda\otimes\sigma\cong\lambda^{(\dim\sigma)}$, so, for any
sequence $\{\sigma_n\}_{n\in\mathbb N}$ of finitedimensional
representations of $\mathbb F_2$, $s\mapsto
\pi_s\otimes(\oplus_{n=1}^\infty\sigma_n)$ is a homotopy that
connects $\lambda^{(\infty)}$ with
$\lambda^{(\infty)}\oplus(\oplus_{n=1}^\infty\sigma_n)$.

It was shown in \cite{Haagerup-T} that there exists a sequence
$\{\sigma_n\}_{n\in\mathbb N}$ of finitedimensional
representations of $\mathbb F_2$ such that
 \begin{equation}\label{Haagerup}
\|\lambda(a)\|=\lim\sup\nolimits_{n\to\infty}\|\sigma_n(a)\|
 \end{equation}
for any $a\in\mathbb C[\mathbb F_2]$. From now on let us fix this
sequence. Put, for $t\in[n,n+1]$,
 $$
\rho_t=(\pi_0\otimes\sigma_1)\oplus\ldots\oplus(\pi_0\otimes\sigma_n)\oplus
(\pi_{t-n}\otimes\sigma_{n+1})\oplus(\pi_1\otimes\sigma_{n+2})\oplus\ldots.
 $$
Then we get a family $(\rho_t)_{t\in[0,\infty)}$ of
representations of $\mathbb F_2$.

Recall that an asymptotic representation of a $C^*$-algebra $A$ is
an asymptotic homomorphism \cite{C-H} into the $C^*$-algebra
$\mathbb B(H)$ of bounded operators on a Hilbert space, i.e. a
family of maps $(\varphi_t)_{t\in[0,\infty)}:A\to\mathbb B(H)$
such that the map $t\mapsto\varphi_t(a)$ is continuous for any
$a\in A$ and $\varphi_t$ behaves asymptotically like a
$*$-homomorphism, as $t\to\infty$. Any continuous family of
(genuine) representations gives an asymptotic representation. A
more general example of an asymptotic representation can be
obtained from a continuous family $\mu_t$, $t\in[0,\infty)$, of
representations of a $C^*$-algebra $E$ such that $q:E\to A$ is a
quotient $*$-homomorphism if this family satisfies the condition
$\lim_{t\to\infty}\|\mu_t(e)\|\leq\|q(e)\|$ for any $e\in E$. If
$\chi:A\to E$ is a Bartle--Graves continuous lifting map \cite{BG}
then $\mu_t\circ\chi:A\to\mathbb B(H)$ is the required asymptotic
representation.

 \begin{thm}\label{Th}
The family $(\rho_t)_{t\in[0,\infty)}$ defines an asymptotic
representation of $C^*_r(\mathbb F_2)$.
 \end{thm}

 \begin{proof}
We are going to show that
$\lim_{t\to\infty}\|\rho_t(a)\|=\|\lambda(a)\|$ for any
$a\in\mathbb C[\mathbb F_2]$. It follows from the definition of
$\rho_t$ that $\lim_{t\to\infty}\|\rho_t(a)\|$ is equal to the
maximum of the following three numbers:
 \begin{eqnarray*}
M_1(a)&=&\lim\nolimits_{n\to\infty}\max\nolimits_{1\leq i\leq
n}\|(\pi_0\otimes\sigma_i)(\Delta (a))\|;\\
M_2(a)&=&\sup\nolimits_{s\in[0,1]}\lim\nolimits_{n\to\infty}\|(\pi_s\otimes\sigma_n)(\Delta(a))\|;\\
M_3(a)&=&\lim\sup\nolimits_{n\to\infty}\|(\pi_1\otimes\sigma_n)(\Delta(a))\|,
 \end{eqnarray*}
where $\Delta$ is the map given by the diagonal homomorphism
$\mathbb F_2\to\mathbb F_2\otimes\mathbb F_2$. It follows from the
Fell's absorbtion principle that
$\pi_0\otimes\sigma_n\cong\lambda^{(\dim\sigma_n)}$ and
$\pi_1\otimes\sigma_n\cong\lambda^{(\dim\sigma_n)}\oplus\sigma_n$,
hence $M_1(a)=\|\lambda(a)\|$ and
 $$
M_3(a)=\max\bigl\lbrace\|\lambda(a)\|,
\,\lim\sup\nolimits_{n\to\infty}\|\sigma_n(a)\|\bigr\rbrace=\|\lambda(a)\|
 $$
by (\ref{Haagerup}), so it remains to estimate $M_2(a)$.

Consider the ideal $I=C^*_{\oplus_{n=1}^\infty\sigma_n}(\mathbb
F_2)\cap\K$ in $C^*_{\oplus_{n=1}^\infty\sigma_n}(\mathbb F_2)$.
The equality (\ref{Haagerup}) implies that
 \begin{equation}\label{isomo}
C^*_r(\mathbb F_2)\cong C^*_{\oplus_{n=1}^\infty\sigma_n}(\mathbb
F_2)/I.
 \end{equation}

By Proposition \ref{exact}, for any $s\in[0,1]$ there is an
$C^*$-algebra isomorphism
 \begin{equation}\label{ravny}
C^*_{\pi_s}(\mathbb
F_2)\otimes(C^*_{\oplus_{n=1}^\infty\sigma_n}(\mathbb F_2)/I)=
(C^*_{\pi_s}(\mathbb F_2)\otimes
C^*_{\oplus_{n=1}^\infty\sigma_n}(\mathbb
F_2))/(C^*_{\pi_s}(\mathbb F_2)\otimes I),
 \end{equation}
where $\otimes$ denotes the minimal tensor $C^*$-product, in other
words, the norm on the left-hand side of (\ref{ravny}) is equal to
the norm on the right-hand side. Let $q$ denote the quotient
$*$-homomorphism from $C^*_{\pi_s}(\mathbb F_2)\otimes
C^*_{\oplus_{n=1}^\infty \sigma_n}(\mathbb F_2)$ to the
$C^*$-algebra in (\ref{ravny}) and let $c\in C^*_{\pi_s}(\mathbb
F_2)\odot C^*_{\oplus_{n=1}^\infty \sigma_n}(\mathbb F_2)$, where
$\odot$ denotes the algebraic tensor product. Then the norm of
$q(c)$ in the right-hand side $C^*$-algebra equals
 $$
\|q(c)\|=\lim\sup\nolimits_{n\to\infty}\|(\pi_s\otimes\sigma_n)(c)\|.
 $$

Now let $c=\Delta(a)$, where $a\in \mathbb C[\mathbb F_2]$ and
$\Delta$ is considered as the diagonal homomorphism,
$\Delta:C^*(\mathbb F_2)\to C^*_{\pi_s}(\mathbb F_2)\otimes
C^*_{\oplus_{n=1}^\infty \sigma_n}(\mathbb F_2)$. Then the norm of
$q(c)=q(\Delta(a))$ in the right-hand side of (\ref{ravny}) equals
 $$
\|q(\Delta(a))\|=\lim\sup\nolimits_{n\to\infty}\|(\pi_s\otimes\sigma_n)(a)\|.
 $$
Because of (\ref{ravny}) we can evaluate the same norm using the
left-hand side of (\ref{ravny}):
 $$
q(\|\Delta(a))\|=\|(\pi_s\otimes\lambda)(a)\|
 $$
(we use here the isomorphism (\ref{isomo})), but the latter equals
$\|\lambda(a)\|$ by the Fell's absorbtion principle, so
 \begin{equation}\label{e1}
\lim\sup\nolimits_{n\to\infty}\|(\pi_s\otimes\sigma_n)(a)\|=
\|\lambda(a)\|.
 \end{equation}
Let us show that convergence in (\ref{e1}) is uniform. This would
follow from equicontinuity of the sequence $\{f_n\}_{n\in\mathbb
N}$ of functions, where $f_n(s)=\|(\pi_s\otimes\sigma_n)(b)\|$,
$s\in[0,1]$, for any finite sum $b=\sum_{i,j}\alpha_{ij}g_i\otimes
g_j\in C^*_{\pi_s}(\mathbb F_2)\odot C^*_{\sigma_n}(\mathbb F_2)$,
where $\alpha_{i,j}\in\mathbb C$, $g_i\in\mathbb F_2$.
Equicontinuity follows from the obvious estimate
 $$
|f_n(s')-f_n(s)|\leq\sum\nolimits_{i,j}|\alpha_{ij}|\cdot\|\pi_{s'}(g_i)-\pi_s(g_i)\|\cdot\|\sigma_n(g_j)\|,
 $$
which is independent of $n$ since $\|\sigma_n(g_j)\|=1$ for all
$n$ and $j$. Finally we can conclude that
 $$
M_2(a)\leq
\sup\nolimits_{s\in[0,1]}\lim\sup\nolimits_{n\to\infty}\|(\pi_s\otimes\sigma_n)(a)\|=\|\lambda(a)\|.
 $$

 \end{proof}

 \begin{dfn}
An asymptotic representation
$(\varphi_t)_{t\in[0,\infty)}:A\to\mathbb B(H)$ is {\it
representation-like} if, for any $a_1,\ldots,a_n\in A$ and for any
$\varepsilon>0$ there exists $T$ such that for any $t>T$ there is
an isometry $U_t$ that satisfies
$\|\varphi_t(a_i)-U_t\pi(a_i)U_t^*\|<\varepsilon$, $i=1,\ldots,n$,
where $\pi:A\to\mathbb B(H)$ is some faithful representation.
 \end{dfn}

We show below that the asymptotic representation
$(\rho_t)_{t\in[0,\infty)}$ from Theorem \ref{Th} is not
representation-like.

\section{Applications}

\subsection{Semi-invertibility of extensions}

Let $E=C^*_{\oplus_{n=1}^\infty\sigma_n}(\mathbb F_2)+\K$. Then
 \begin{equation}\label{extension}
0\to\K\to E\to E/\K=C^*_{\oplus_{n=1}^\infty\sigma_n}(\mathbb
F_2)/I\cong C^*_r(\mathbb F_2)\to 0
 \end{equation}
is an extension of $C^*_r(\mathbb F_2)$ by $\K$ constructed in
\cite{Haagerup-T}, where it was proved that this extension is
non-invertible. Here we show that this extension is
semi-invertible \cite{MT3}, moreover, it represents the trivial
element in the group $\Ext^{-1/2}(C^*_r(\mathbb F_2),\K)$
\cite{MT3}. Let us briefly recall these definitions. An extension
$0\to \K\to E\stackrel{p}\to A\to 0$, of $A$ by $\K$, is {\it
asymptotically split} if there exists an asymptotic homomorphism
$(\varphi_t)_{t\in[0,\infty)}:A\to E$ such that
$p\circ\varphi_t=\id_A$ for any $t\in[0,\infty)$. An extension of
$A$ by $\K$ is {\it semi-invertible} if its direct sum with some
other extension of $A$ by $\K$ is asymptotically split. Two
extensions of $A$ by $\K$ are {\it stably equivalent} if their
direct sums with some semi-invertible extension of $A$ by $\K$ are
unitarily equivalent. Finally, the group $\Ext^{-1/2}(A,\K)$ is
the group of classes of stable equivalence of semi-invertible
extensions of $A$ by $\K$. The only difference of
$\Ext^{-1/2}(A,\K)$ from the group $\Ext(A,\K)^{-1}$ of
invertibles in the classical BDF functor is that we use
asymptotically split extensions instead of split ones.

 \begin{cor}
The extension $($\ref{extension}$)$ is semi-invertible. Moreover,
it represents the trivial element in the group
$\Ext^{-1/2}(C^*_r(\mathbb F_2),\K)$.
 \end{cor}
 \begin{proof}
Let $q:\mathbb B(H)\to \mathbb Q=\mathbb B(H)/\K$ be the quotient
map. Then an easy calculation shows that
 $$
q(\rho_t(a))=q\bigl(\oplus_{n=1}^\infty\sigma_n(a)\oplus\lambda^{(\infty)}(a)\bigr).
 $$
Therefore, the direct sum of the extension (\ref{extension}) and
of the trivial extension admits an asymptotic lift.
 \end{proof}

This gives an evidence to the following conjecture:
 \begin{conj}
All extensions of $C^*_r(\mathbb F_2)$ are semi-invertible.
 \end{conj}


\subsection{Estimates for the asymptotic tensor $C^*$-norm}

In \cite{MT6} we have introduced the asymptotic tensor $C^*$-norm
(one-sided and symmetric) for tensor products of separable
$C^*$-algebras. Recall that the left asymptotic tensor $C^*$-norm
on the algebraic tensor product $A\odot B$ of two $C^*$-algebras
is defined by
 $$
\|c\|_\lambda=\sup\nolimits_{\varphi,\psi}\lim\sup\nolimits_{t\to\infty}\|(\varphi_t\otimes\psi)(c)\|,
 $$
where $c\in A\odot B$ and the supremum is taken over all
asymptotic representations
$\varphi=(\varphi_t)_{t\in[0,\infty)}:A\to\mathbb B(H)$ and all
representations $\psi:B\to\mathbb B(H)$. The symmetric asymptotic
tensor norm $\|\cdot\|_\sigma$ is defined as a similar supremum,
where $\psi$ stands for asymptotic representations of $B$.

Now we are able to calculate the left asymptotic norm of some
elements in $C^*_r(\mathbb F_2)\odot
C^*_{\oplus_{n=1}^\infty\bar{\sigma}_n}(\mathbb F_2)$, where
$\bar{\sigma}$ is the representation contragredient to $\sigma$,
and the symmetric asymptotic tensor $C^*$-norm of some elements in
$C^*_r(\mathbb F_2)\odot C^*_r(\mathbb F_2)$.

Let $a=\frac{1}{4}(g_1+g_1^{-1}+g_2+g_2^{-1})$, where $g_1,g_2$
are free generators for $\mathbb F_2$.

 \begin{prop}\label{estim}
$\|\Delta(a)\|_\lambda=1$ in $C^*_r(\mathbb F_2)\odot
C^*_{\oplus_{n=1}^\infty\bar{\sigma}_n}(\mathbb F_2)$;
$\|\Delta(a)\|_\sigma=1$ in $C^*_r(\mathbb F_2)\odot C^*_r(\mathbb
F_2)$.
 \end{prop}
 \begin{proof}
Obviously $\|\Delta(a)\|_\lambda\leq \|a\|\leq 1$. Similarly,
$\|\Delta(a)\|_\sigma\leq 1$. Let us use the asymptotic
representation $\rho_t$ to estimate both norms from below.
 $$
\|\Delta(a)\|_\lambda\geq\lim\sup\nolimits_{t\to\infty}
\|(\rho_t\otimes(\oplus_{n=1}^\infty\bar{\sigma}_n))(\Delta(a))\|\geq
\lim\sup\nolimits_{n\to\infty}\|(\sigma_n\otimes\bar{\sigma}_n)(\Delta(a))\|=1
 $$
because $\sigma_n\otimes\bar{\sigma}_n$ contains a copy of the
trivial representation.

Similarly,
 $$
\|\Delta(a)\|_\sigma\geq\lim\sup\nolimits_{t\to\infty}
\|(\rho_t\otimes\bar{\rho}_t)(\Delta(a))\|\geq
\lim\sup\nolimits_{n\to\infty}\|(\sigma_n\otimes\bar{\sigma}_n)(\Delta(a))\|=1.
 $$
 \end{proof}

 \begin{cor}
The asymptotic representation of Theorem \ref{Th} is not
representation-like.
 \end{cor}
 \begin{proof}
It is easy to see that for a representation-like asymptotic
representation $(\varphi_t)_{t\in[0,\infty)}:A\to\mathbb B(H)$ one
has
 $$
\lim\sup\nolimits_{t\to\infty}\|(\varphi_t\otimes\mu)(c)\|\leq
\|(\pi\otimes\mu)(c)\|
 $$
for any $C^*$-algebra $B$, any representation $\mu$ of $B$ and any
$c\in A\odot B$, where $\pi$ is (any) faithful representation of
$A$.

 \end{proof}

The following corollary shows that there is no asymptotic analogue
for the Fell's absorbtion principle.
 \begin{cor}
$C^*_r(\mathbb F_2)\otimes_\sigma C^*_r(\mathbb F_2)\neq
C^*_r(\mathbb F_2)\otimes C^*_r(\mathbb F_2)$.
 \end{cor}
 \begin{proof}
It is well-known \cite{Kesten} that
$\|(\lambda\otimes\lambda)(\Delta(a))\|=\|\lambda(a)\|=\frac{\sqrt{3}}{2}<1=\|\Delta(a)\|_\sigma$.
 \end{proof}

Note that $C^*_r(\mathbb F_2)\otimes_\lambda C^*_r(\mathbb
F_2)\neq C^*_r(\mathbb F_2)\otimes C^*_r(\mathbb F_2)$ due to the
Fell's absorbtion principle.

Note also that the result of Proposition \ref{estim} contrasts the
analogous result for property T groups, cf. \cite{MT7}. Let $G$ be
a property T group, $g_1,\ldots,g_m$ its generators and
$a=\frac{1}{2m}(g_1+g_1^{-1}+\ldots+g_m+g_m^{-1})\in\mathbb C[G]$,
$\alpha>0$ the Kazhdan constant (i.e.
$\Sp(a)\subset[-1,1-\alpha]\cup\{1\}$, cf. \cite{HRV}). Let
$\tau_n$, $n\in\mathbb N$, be a sequence of finitedimensional
representations of $G$. The proof of the following proposition is
essentially contained in \cite{MT7}.

 \begin{prop}\label{estim2}
$\|\Delta(a)\|_\lambda=1-\alpha$ in $C^*_r(G)\odot
C^*_{\oplus_{n=1}^\infty\tau_n}(G)$.
 \end{prop}
 \begin{proof}
It is known \cite{Valette} that $\Delta(a)=x+\Delta(p)$, where
$x\leq 1-\alpha$ and $p\in C^*(G)$ is the projection corresponding
to the trivial representation of $G$. Let
$(\varphi_t)_{t\in[0,\infty)}:C^*_r(G)\to\mathbb B(H)$ be an
arbitrary asymptotic representation. The proposition would follow
if we show that
 \begin{equation}\label{nol}
\lim_{t\to\infty}\bigl(\varphi_t\otimes(\oplus_{n=1}^\infty\tau_n)\bigr)(\Delta(p))=0.
 \end{equation}

Put
 $$
q_t(n)=(\varphi_t\otimes\tau_n)(\Delta(p));\quad
q_t=\oplus_{n=1}^\infty q_t(n).
 $$
Then $q_t$ is asymptotically a projection, i.e.
$\lim_{t\to\infty}\|q_t-q_t^2\|=0$, hence the same holds for each
$n$ uniformly: $\lim_{t\to\infty}\|q_t(n)-q_t^2(n)\|=0$. This
means that one can determine if $\lim_{t\to\infty}q_t(n)=0$ by
looking at $q_t(n)$ for some fixed $t$ (same for all $n$).

Suppose that (\ref{nol}) is not true. Then there exists $n$ such
that $\lim_{t\to\infty}\|q_t(n)\|=1$. Hence
$\lim_{t\to\infty}\|(\varphi_t\otimes\tau_n)(\Delta(a))\|=1$, i.e.
the norm of $\Delta(a)$ in $C^*_r(G)\otimes_\lambda
C^*_{\tau_n}(G)$ equals 1. But since $C^*_{\tau_n}(G)$ is
finitedimensional, we have $C^*_r(G)\otimes_\lambda
C^*_{\tau_n}(G)=C^*_r(G)\otimes C^*_{\tau_n}(G)$. But the norm of
$\Delta(a)$ in $C^*_r(G)\otimes C^*_{\tau_n}(G)$ cannot equal 1
\cite{WassT}. This contradiction ends the proof.

 \end{proof}

\end{document}